\documentclass[12pt]{amsart}
\usepackage{amscd,amssymb}
\usepackage[graph,frame,poly,arc]{xy}  
\usepackage[plainpages,backref,urlcolor=blue]{hyperref}
\usepackage{mathtools}

\theoremstyle{plain}
\newtheorem{thm}[subsection]{Theorem}

\newtheorem{prop}[subsection]{Proposition}
\newtheorem{cor}[subsection]{Corollary}

\theoremstyle{definition}
\newtheorem{rk}[subsection]{Remark}
\newtheorem{definition}[subsection]{Definition}
\newtheorem{ex}[subsection]{Example}

\newtheorem{question}[subsection]{Question}

\numberwithin{equation}{section}
\newcommand{\CC}{{\mathcal C}}

\newcommand{\al}{{\alpha}}
\newcommand{\be}{{\beta}}

\newcommand{\C}{\mathbb{C}}
\newcommand{\PP}{\mathbb{P}}

\DeclareMathOperator{\mult}{mult}

\begin{document}

\title [Plus-one generated curves, Brian\c con polynomials and eigenschemes]
{Plus-one generated curves, Brian\c con-type polynomials and eigenscheme ideals}

\author[Alexandru Dimca]{Alexandru Dimca}
\address{Universit\'e C\^ ote d'Azur, CNRS, LJAD, France and Simion Stoilow Institute of Mathematics,
P.O. Box 1-764, RO-014700 Bucharest, Romania}
\email{dimca@unice.fr}

\author[Gabriel Sticlaru]{Gabriel Sticlaru}
\address{Faculty of Mathematics and Informatics,
Ovidius University
Bd. Mamaia 124, 900527 Constanta,
Romania}
\email{gabriel.sticlaru@gmail.com }



\begin{abstract} We define the minimal plus-one generated curves and prove a result explaining why they are the closest relatives of the free curves, after the nearly free curves. Then we look at the projective closures of the general and of the special fibers of some Briançon-type polynomials constructed by E. Artal Bartolo, Pi. Cassou-Noguès and  I. Luengo Velasco. They yield new examples of free, nearly free or
minimal plus-one generated curves, as well as counter-examples to the conjecture saying that a supersolvable curve is free.
In the final section we give a characterization of plus-one generated curves in terms of eigenscheme ideals, similar to the characterization of free curves given by 
R. Di Gennaro, G. Ilardi, R.M. Miró-Roig, H. Schenck and J. Vallès in a recent paper. Then we apply this result to the construction of minimal plus-one generated curves obtained by putting together at least two members in a pencil of curves related to Briançon-type polynomials.
\end{abstract}
 
\maketitle


\section{Introduction}

Let $S=\C[x,y,z]$ be the polynomial ring in three variables $x,y,z$ with complex coefficients, and let $C:f=0$ be a reduced curve of degree $d\geq 3$ in the complex projective plane $\PP^2$. 
We denote by $J_f$ the Jacobian ideal of $f$, i.e., the homogeneous ideal in $S$ spanned by the partial derivatives $f_x,f_y,f_z$ of $f$, and  by $M(f)=S/J_f$ the corresponding graded quotient ring, called the Jacobian (or Milnor) algebra of $f$.
The $S$-module of derivations $Der(S)$ of the polynomial ring $S$ is a free module with a basis given by $\partial_x$, $\partial_y$ and $\partial_z$, that is any derivation $\theta \in Der(S)$ can be uniquely written as
  \begin{equation}
\label{eqE0}
 \theta=a\partial_x +b \partial_y+ c \partial_z
\end{equation} 
 with $a,b,c \in S$. When  $a,b,c \in S_m$, we say that $\theta$ has degree $m$ and write $\deg \theta=m$. 
Consider the graded $S$-module of Jacobian syzygies of $f$ or, equivalently, the module of derivations killing $f$, namely
\begin{equation}
\label{eqD0}
D_0(f)= \{\theta \in Der(S) \ : \ \theta(f)=0\}= \{\rho=(a,b,c) \in S^3 \ : \ af_x+bf_y+cf_z=0\}.
\end{equation}
 The  graded $S$-module
 $$D(f)=\{\theta \in Der(S) \ : \ \theta (f) \in (f) \}$$
 of derivations preserving $f$ has a direct sum decomposition in each degree $s$ given by
 \begin{equation}
\label{eqE1}
D(f)_s=S_{s-1}\cdot E \oplus D_0(f)_s,
\end{equation}
for any integer $s$, where $E= x\partial_x +y\partial_y+ z\partial_z$ is the Euler derivation and $D_0(f)$ is the submodule of derivations killing $f$ as in \eqref{eqD0}.
We say that $C:f=0$ is an {\it $m$-syzygy curve} if  the module $D_0(f)$ is minimally generated by $m$ homogeneous syzygies, say $\rho_1,\rho_2,\ldots ,\rho_m$, of degrees $d_j=\deg \rho_j$ ordered such that $$d_1\leq d_2 \leq \ldots \leq d_m.$$ 
In this note we assume  that $C$ is not the union of $d$ lines passing through one point, hence $d_1 \geq 1$.
We call these degrees $(d_1, \ldots, d_m)$ the {\it exponents} of the curve $C$. 
The smallest degree $d_1$ is sometimes denoted by ${\rm mdr}(f)$ and is called the {\it minimal degree of a Jacobian relation} for $f$.

The curve $C$ is {\it free} when $m=2$, since then  $D_0(f)$ is a free module of rank 2, see for instance \cite{KS,Sim2,ST,To}.
A free curve $C$ is also characterized by the condition $d_1+d_2=d-1$, and hence in this case
$r=d_1 <d/2$, see \cite{ST}.  A {\it nearly free} curve, as introduced in \cite{DStRIMS}, is a 3-syzygy curve $C$ with exponents $(d_1,d_2,d_3)$ such that $d_1+d_2=d$ and
$d_2=d_3$. To simplify the notation, we say that $(d_1,d_2)$ are the exponents for the nearly free curve $C$.
Finally, a plus-one generated curve is  a 3-syzygy curve $C$ with exponents $(d_1,d_2,d_3)$ such that $d_1+d_2=d$, see
\cite{Abe,minTjurina}. We define the $\delta$-{\it level} of a  plus-one generated curve by the formula
\begin{equation}
\label{eqD1}
\delta L(C)=d_3-d_2 \geq 0.
\end{equation}
Note that a plus-one generated curve $C$ is nearly free if and only if $\delta L(C)=0$. 
\begin{definition}
\label{def1}
A plus-one generated curves $C$ satisfying $\delta L(C)=1$ is called  a {\it minimal plus-one generated } curve, or an $MPOG$-curve for short.
\end{definition}
These curves can be regarded as the closest neighbors of the free curves beyond the nearly free curves. Before stating the result, we recall some notions and results.
The total Tjurina number $\tau(C)$ of a reduced plane curve $C:f=0$ is just the degree of its Jacobian ideal $J_f$, or equivalently the sum of the Tjurina numbers of all the singularities of $C$.
The following result,  due to  du Plessis and Wall, see \cite[Theorem 3.2]{duPCTC} as well as \cite{E} for an alternative approach, gives bounds on this Tjurina number $\tau(C)$.
\begin{thm}
\label{thmCTC}
For positive integers $d$ and $r$, define two new integers by 
$$\tau(d,r)_{min}=(d-1)(d-r-1)  \text{ and } 
\tau(d,r)_{max}= (d-1)^2-r(d-r-1).$$ 
Then, if $C:f=0$ is a reduced curve of degree $d$ in $\PP^2$ and  $r={\rm mdr}(f)$,  one has
$$\tau(d,r)_{min} \leq \tau(C) \leq \tau(d,r)_{max}.$$
Moreover, for $r={\rm mdr}(f) \geq d/2$, the stronger inequality
$\tau(C) \leq \tau'(d,r)_{max}$ holds, where
$$\tau(d,r)'_{max}=\tau(d,r)_{max} - \binom{2r+2-d}{2}.$$

\end{thm}
A degree $d$ curve $C:f=0$ satisfying $r={\rm mdr}(f) \geq d/2$ and 
$\tau(C)= \tau'(d,r)_{max}$ is called a {\it maximal Tjurina curve of type} $(d,r)$, see \cite{maxTjurina}.

At the end of the proof of Theorem \ref{thmCTC}, in \cite{duPCTC}, the authors state the following very interesting consequence (of the proof, not of the statement) of Theorem \ref{thmCTC}.
\begin{cor}
\label{corCTC} Let $C:f=0$ be a reduced curve of degree $d$ in $\PP^2$ and  $r={\rm mdr}(f)$. One has
$$ \tau(C) =\tau(d,r)_{max}$$
if and only if $C:f=0$ is a free curve, and then $r <d/2$.
\end{cor}
In the paper \cite{Dmax}, the first author has given an alternative proof of Corollary \ref{corCTC} and has
shown that the following similar property holds for nearly free curves.

\begin{cor}
\label{corCTC2} Let $C:f=0$ be a reduced curve of degree $d$ in $\PP^2$ and  $r={\rm mdr}(f)$.
One has
$$ \tau(C) =\tau(d,r)_{max}-1$$
if and only if $C:f=0$ is a nearly free curve, and then $r  \leq d/2$.
\end{cor}
Our first main result in this paper is the following.
\begin{thm}
\label{thm0} Let $C:f=0$ be a reduced curve of degree $d$ in $\PP^2$ and  $r={\rm mdr}(f)$.
Then  $C$ is an $MPOG$-curve if and only if
$$ \tau(C) =\tau(d,r)_{max}-2.$$
\end{thm}
Let $I_f$ denote the saturation of the ideal $J_f$ with respect to the maximal ideal ${\bf m}=(x,y,z)$ in $S$ and consider  the graded {\it Jacobian module} of $f$ defined by 
 $$N(f)=I_f/J_f.$$
We set $n(f)_k=\dim N(f)_k$ for any integer $k$ and introduce the {\it freeness defect of the curve} $C$ by the formula
$$\nu(C)=\max _j \{n(f)_j\}$$ as in \cite{AD}.
\begin{rk}
\label{rkNU}
Note that $C$ is free if and only if $N(f)=0$, see for instance \cite{ST}, and hence in this case $\nu(C)=0$, and $C$ is nearly free if and only if $\nu(C)=1$, see \cite{DStRIMS}. To prove Theorem \ref{thm0}, we show in fact that $ \tau(C) =\tau(d,r)_{max}-2$ implies $\nu(C)=2$ and then use the classification of such curves given in \cite[Theorem 3.11]{minTjurina}. This result says that a curve $C$ with $\nu(C)=2$ is either an MPOG-curve, or a 4-syzygy curve with exponents $(r,r,r,r)$ and degree
$d=2r-1$, which is in fact a {\it maximal Tjurina curve of type} $(d,r)=(2r-1,r)$
as defined above.
\end{rk}

The {\it Brian\c con-type polynomials} are non-homogeneous polynomials $h$
in the ring of polynomials
$R=\C[x,y]$, such that the associated mapping $h: \C^2 \to \C$ has all the fibers $h^{-1}(t)$ smooth and irreducible, and yet $h$ is not obtained from a linear form in $R$ by composition with a polynomial automorphism of $\C^2$. A series of such polynomials have been constructed in \cite{ACL}, the simplest of them being two polynomials of degree 10 that we recall now.

Let  $s=xy+1$, $p=xs+1$ and $u=s^2+y$. Using this notation, we define 
\begin{equation}
\label{eqB1}
g=p^2u-\frac{5}{3}ps-\frac{1}{3}s \text{ and } 
g'=p^2u-\frac{7}{9}ps+\frac{1}{9}s.
\end{equation}
The polynomial $g$ up-to a change of coordinates $(x,y) \mapsto (-x,-y)$ was considered in an unpublished manuscript by J. Brian\c con, see \cite{ACL, ACD}.
We recall that for a polynomial mapping $h:\C^2 \to \C$ the {\it set of atypical values} is the minimal set
$A(h) \subset \C$ such that $h$ induces a locally trivial topological fibration over the set $\C \setminus (A(h) \cup C(h))$, where $C(h)$ denotes the set of critical values of $h$. The atypical values for a function $h:\C^2 \to \C$ come from jumping of the Milnor numbers of the singularities of $h$ at infinity, see \cite{Bt, Dur, HL}.
The polynomials $g$ and $g'$ have the following very interesting properties, see \cite{ACL}.
 \begin{thm}
\label{thmB}

\begin{enumerate}

\item The polynomials $g$ and $g'$ have no critical points, that is all the fibers of $g$ and $g'$ are smooth affine curves. In other words, the sets $C(g)$ and $C(g')$ are empty.

\item All the fibers of the polynomials $g$ and $g'$ are irreducible affine plane curves.

\item The set of atypical values for $g$ is 
$A(g)= \{ 0, b \}$, where $b= -\frac{16}{9}.$

\item The set of atypical values for $g'$ is 
$A(g')= \{ 0, b'  \}$, where $b'=-\frac{64}{81}.$

\item The generic fiber $F_1=g^{-1}(1)$ of the polynomial $g$ and 
$F_1'=g'^{-1}(1)$  of the polynomial $g$ have their first Betti numbers given by $b_1(F_1)=b_1(F_1')=4$. In particular, the polynomials $g$ and $g'$ are not obtained from linear forms by composition with a polynomial automorphism of $\C^2$.

\end{enumerate}

\end{thm}

When $h: \C^2 \to \C$ is obtained from a linear form in $R$ by composition with a polynomial automorphism of $\C^2$, the freeness properties of the projective closure $D_t$ of any fiber $h^{-1}(t)$ are discussed
in \cite[Theorem 1.5]{POG}. It is natural to ask what happens when we replace such a polynomial $h$ by a Brian\c con-type polynomial.

For the polynomial $g$, we consider three associated plane curves $C_{t}$ in $\PP^2$, which are the projective closure of the affine fibers $F_t=g^{-1}(t)$ for $t=0$, $t=b$ and $t=1$.
The first two fibers correspond to the atypical values in $A(g)$, and the third one is a general fiber of $g$. For the convenience of the reader, the defining equations for these curves, obtained by homogenization starting with the first equation in \eqref{eqB1}, are given in \eqref{eqB2} and
 \eqref{eqB3}.
Similarly, for the polynomial $g'$, we consider three associated plane curves $C_t'$ in $\PP^2$, which are the projective closure of the affine fibers $g'^{-1}(t)$ for $t=0$, $t=b'$ and $t=1$. The corresponding equations are given in \eqref{eqNB1} and
 \eqref{eqNB2}.

 Our results about the freeness properties of these projective plane curves are given in Proposition \ref{propthm1} and Proposition  \ref{propthm2}.

 Since the line at infinity $L_z:z=0$ intersects the curves $C_{t}$ only in two points, the general principle explained in \cite{DIPS} tells us that the curves $C_{t} \cup L_z$  (resp. $C_{t} '\cup L_z$)  of degree $11$ may also enjoy (sometimes better) freeness properties. See also Step 3 in the proof of Theorem 1.5 in \cite{POG} for such a situation. We consider also the curve arrangements  $C_0 \cup C_{b}$ (resp. $C_0' \cup C'_{b'}$)  and $C_0 \cup C_{b}\cup L_z$ (resp. $C_0' \cup C'_{b'}\cup L_z$.
The corresponding results are given in Propositions \ref{propthm3},  \ref{propthm4} and  \ref{propthm5}.

At the end of Section 3, we construct two {\it counter-examples} involving the curve $C_0$ to the conjecture \cite[Conjecture 1.10]{DIPS}, which asks whether {\it any supersolvable curve is free}. In the first example
the modular point is on $C_0$, while in the second example the modular point is not on $C_0$, see Remark \ref{rkB10}, where we recall
the definitions of all the notions involved here.

 Using the results in \cite{ACL}, one can construct other new curves of higher degree, which are either free, or nearly free or  MPOG-curves, see Section 4 below. This leads us to ask the natural question if such a result holds for any Brian\c con-type polynomials.
 
   \begin{question}
\label{Q1} 
Do the projective closures of the general and the special (atypical) fibers of any
Brian\c con-type polynomial enjoy some kind of freeness property, i.e., do the corresponding exponents have some special properties ?
 \end{question}
It seems that this question is related to the very challenging question of finding numerical
restrictions for Brian\c con-type polynomials, involving their geometry.

Finally we give a characterization of plus-one generated curves in terms of eigenscheme ideals, similar to the following characterization of free curves given by 
R. Di Gennaro, G. Ilardi, R.M. Miró-Roig, H. Schenck and J. Vallès in a recent paper \cite{Di+}.
For a derivation $\theta$ as in \eqref{eqE0}, we consider the ideal
$I_{\theta} \subset S$ generated by the three 2-minors of the matrix having as a first row $(x,y,z)$ and as a second row the coefficients $(a,b,c)$ of $\theta$. The {\it eigenscheme} $\Gamma _{\theta}$ associated to $\theta$ is the subscheme in $\PP^2$ defined by the ideal $I_{\theta}$. When  the eigenscheme $\Gamma _{\theta}$ is 0-dimensional, then the ideal $I_{\theta}$ is {\it saturated}, see \cite[Proposition 3.1]{Be}.
The following result was obtained in \cite{Di+}, see Theorem 2.5.
 \begin{thm}
\label{thmE1}
Let $C:f=0$ be a reduced curve in $\PP^2$ of degree $d$. Let $\theta \in D_0(f)$ be a non trivial derivation of degree $r$. We assume that $2r \leq d-1$ and that the eigenscheme $\Gamma _{\theta}$ is 0-dimensional.
Then $C$ is a free curve with exponents $(r,d-r-1)$ if and only if $f \in I_{\theta}$.
 \end{thm}
 The condition $2r \leq d-1$ is added here since we like our exponents $(d_1,d_2)$ to verify the property $d_1 \leq d_2$. See also Theorem \ref{thmE2} $(1)$ below. 
In view of the direct sum decomposition, one may even take $\theta \in D(f)_r$, as in \cite[Theorem 2.5]{Di+}, and the same result holds.
A key remark is that $f \in I_{\theta}$ if and only if the ideal quotient
$$K_{\theta,f}=I_{\theta}:(f)=\{h \in S \ : \ hf \in I_{\theta} \} $$
is the whole ring $S$,
where $(f)$ denotes the principal ideal in $S$ generated by $f$. From this new point of view, Theorem \ref{thmE1} can be modified to cover the case of plus-one generated curves, as our second main result shows.
 \begin{thm}
\label{thmE2}
Let $C:f=0$ be a reduced curve in $\PP^2$ of degree $d$. Let $\theta \in D_0(f)$ be a non trivial derivation of degree $r$. We assume that $2r \leq d-1$ and that the eigenscheme $\Gamma _{\theta}$ is 0-dimensional.
Then the following hold.
\begin{enumerate}

 \item $r={\rm mdr}(f)$.

\item $C$ is a plus-one generated curve  if and only if the ideal quotient
$K_{\theta,f}=I_{\theta}:(f)$
is a proper ideal having the degree one  component $(K_{\theta,f})_1$ nonzero. 

\end{enumerate}
If the properties in $(2)$ hold, then the ideal $K_{\theta,f}$
has two generators, a linear form $\ell \in S_1$ and a homogeneous
polynomial $h \in S_e$, with $e \geq 1$, and the exponents of the curve $C$ are $(r,d-r,d-r-1+e)$. In particular,
$C$ is nearly free if and only if $e=1$ and $C$ is an MPOG-curve if and only if $e=2$.
 \end{thm}

This result is applied to the construction of minimal plus-one generated curves obtained by putting together some members in a pencil of curves, see Corollary \ref{corE2} for the general situation, and Proposition \ref{propthmE3}, for an explicit construction related to Brian\c con-type polynomials. To our knowledge, this is the first example where putting together the most singular members in a pencil of plane curves does not produce a free curve, see also \cite[Remark 3.1]{Di+}.

\medskip

The first author thanks  Pierrette Cassou-Nogu\` es for useful discussion on this paper, reflected in part in Remark \ref{rkE2a}.

\section{The proof of Theorem \ref{thm0}}

If we assume that  $C$ is an $MPOG$-curve, then the equality for $\tau(C)$ follows from \cite[Proposition 2.1, (4)]{minTjurina}.

Conversely, assume now that $ \tau(C) =\tau(d,r)_{max}-2.$ 
If $r >d/2$, then Theorem \ref{thmCTC} implies that
$$\tau(C) \leq \tau(d,r)_{max}' \leq \tau(d,r)_{max}-\binom{3}{2}=\tau(d,r)_{max}-3,$$
a contradiction. It follows that $r \leq d/2$. Now \cite[Theorem 1.2]{Drcc} implies that 
$$\nu(C)= \tau(d,r)_{max} -\tau(C)=2,$$
when $r<d/2$. In the case $r=d/2$, we have the following equalities, as in the discussion after
\cite[Theorem 1.2]{Drcc}.
$$\nu(C)=3r^2-3r+1- \tau(C)=2,$$
since $\tau(C)=\tau(d,r)_{max}-2=3r^3-3r-1.$
Therefore the assumption $\tau(C)=\tau(d,r)_{max}-2$ implies in all cases that $\nu(C)=2$. Now we apply \cite[Theorem 3.11]{minTjurina} and conclude that $C$ is an $MPOG$-curve.


\begin{rk}
\label{rkG0}
Let $C:f=0$ be a reduced curve of degree $d$ in $\PP^2$ and  $r={\rm mdr}(f)$.
Then the equality
$$ \tau(C) =\tau(d,r)_{max}-3$$
does not determine uniquely the exponents of $C$. Indeed, such a curve may be a plus-one generated curve with exponents $(r,d-r,d-r+2)$, a 4-syzygy curve with exponents $(r,r,r,r)$ and degree $d=2r-1$ as in \cite[Theorem 3.11]{minTjurina}, or maybe even other types of curves.
 
 \end{rk}

\begin{rk}
\label{rkG}
 
 The free and nearly free curves have strong relations with the rational cuspidal plane curves, see \cite{DStRIMS,Mosk}. The  $MPOG$-curves  are related with the rational nearly cuspidal curves, i.e. rational curves whose singularities are all irreducible except one singularity, which has two branches, see \cite[Theorem 1.3, Conjecture 1.4, Corollary 5.8]{minTjurina}.
 See also \cite{MaVa} for an alternative view on nearly free curves.
 \end{rk}


\section{On the freeness of the curves associated to Brian\c con-type polynomials of degree 10}

Here are the defining equations for the curves $C_t$ in $\PP^2$, obtained by homogenization starting with the first equation in \eqref{eqB1}.

\begin{equation}
\label{eqB2}
C_{0}:f_{0}=x^6y^4+4x^5y^3z^2+3x^4y^3z^3+6x^4y^2z^4+\frac{19}{3}x^3y^2z^5+4x^3yz^6+3x^2y^2z^6+\end{equation}
$$+\frac{11}{3}x^2yz^7+x^2z^8+2xyz^8+\frac{1}{3}xz^9+yz^9-z^{10}=0,$$

\begin{equation}
\label{eqB3}
C_{b}:f_{b}= f_{0} +\frac{16}{9}z^{10}=0 \text{  and  } C_{1}:f_{1}=f_{0}-z^{10}=0.
\end{equation}

\begin{rk}
\label{rkB1}
It is clear that all the curves $C_{t}$ have two points situated on the line at infinity $L:z=0$, namely the points $p=(1:0:0)$ and $q=(0:1:0)$.
It is easy to check using SINGULAR \cite{Sing} that one has
$$\mu(C_{0},q)=\mu(C_{b},q)=\mu(C_{1},q)=42$$
and
$$\tau(C_{0},q)=\tau(C_{b},q)=\tau(C_{1},q)=35.$$
Hence at the point $q$ the family of plane curve singularities $(C_{t},q)$ is 
$\mu$-constant and $\tau$-constant. Moreover, all these curve singularities have just one branch, i.e., they are irreducible.
On the other hand we have
$$\mu(C_{0},p)= 28, \ \mu(C_b,p)=30 \text{ and }\mu(C_{1},p)=27.$$
The Tjurina numbers also vary, as we have
$$\tau(C_{0},p)=\tau(C_{1},p)=24 \text{ and } \tau(C_b,p)=26.$$
{\it It is interesting that the total Tjurina number of the exceptional
curve $C_0$ coincides with the total Tjurina number of the generic curve
$C_1$. } Indeed, recall that the exceptional fibers $F_t$ of a polynomial mapping $g:\C^2 \to \C$ are detected by the jump in the total Milnor number of the compactifications $C_t$, see \cite{HL}. 
The number of branches $r(C_t,p)$ for the singularity $(C_t,p)$ is given by the following
$$r(C_0,p)=3, \  r(C_b,p)=1 \text{ and } r(C_1,1)=2.$$
In fact, the claims about the Milnor numbers $\mu(C_t,p)$ and $\mu(C_t,q)$, and about the number of branches $r(C_t,p)$ and $r(C_t,q)$ are also proven in a different way in \cite{ACL}.
\end{rk}

Here are the defining equations for the curves $C'_t$, obtained by homogenization starting with the second equation in \eqref{eqB1}.

\begin{equation}
\label{eqNB1}
C_0':f_{0}'= x^6y^4+4x^5y^3z^2+3x^4y^3z^3+6x^4y^2z^4+\frac{65}{9}x^3y^2z^5+4x^3yz^6+3x^2y^2z^6+
\end{equation}
$$+\frac{49}{9}x^2yz^7+x^2z^8+\frac{10}{3}xyz^8+\frac{11}{9}xz^9+yz^9+\frac{1}{3}z^{10} =0,$$

\begin{equation}
\label{eqNB2}
C_{b'}':f_{b'}'= f_{0}' +\frac{64}{81}z^{10}=0 \text{  and  } C_{1}':f_{1}'= f_{0}' -z^{10}=0.
\end{equation}

\begin{rk}
\label{rkNB1}
It is clear that all the curves $C_t'$ have two points situated on the line at infinity $L:z=0$, namely the points $p=(1:0:0)$ and $q=(0:1:0)$.
It is easy to check using SINGULAR that one has
$$\mu(C_{0}',q)=\mu(C_{b'}',q)=\mu(C_1',q)=42$$
and
$$\tau(C_{0}',q)=\tau(C_{b'}',q)=\tau(C_{1}',q)=35.$$
Hence at the point $q$ the family of plane curve singularities $(C_{t}',q)$ is 
$\mu$-constant and $\tau$-constant. Moreover, all these curve singularities have just one branch.
On the other hand we have
$$\mu(C_0',p)=  \mu(C_{b'}',p)=29 \text{ and }\mu(C_1',p)=27.$$
The Tjurina numbers also vary, as we have
$$\tau(C_0',p)= \tau(C_{b'}',p)=25 \text{ and }\tau(C_1',p)=24.$$
Moreover, all these curve singularities at $p$ have two branches.
These results, except the computation of the Tjurina numbers, are also obtained in an alternative way in \cite{ACL}.
\end{rk}

 \begin{prop}
\label{propthm1}
The geometric genera of the curves $C_t$ are given by
$$g(C_1)=1 \text{ and } g(C_0)=g(C_b)=0.$$
Moreover, these curves of degree $10$  enjoy the following freeness properties.

\begin{enumerate}

\item The curves $C_0$ and $C_1$ are $MPOG$-curves with exponents $(5,5,6)$ and $\tau(C_0)= \tau(C_1)= 59$. 

\item The curve $C_b$ is a free curve with exponents $(4,5)$ and $\tau(C_b)=61$.

\end{enumerate}

 \end{prop}

 \begin{prop}
\label{propthm2}

The geometric genera of the curves $C_t'$ are given by
$$g(C_1')=1 \text{ and } g(C_0')=g(C_{b'}')=0.$$
Moreover, these curves of degree $10$  enjoy the following freeness properties.

\begin{enumerate}

\item The curve  $C_1'$ is an $MPOG$-curve with exponents $(5,5,6)$ and $\tau(C_1')= 59$.

\item The curves $C_0'$ and  $C_{b'}'$ are nearly free curves with exponents $(5,5)$ and $\tau(C_0')=\tau(C_{b'}')=60$.

\end{enumerate}

 \end{prop}
 
 \begin{rk}
\label{rkCC} The claims in Proposition \ref{propthm1} and Proposition \ref{propthm2}
about the genera of the curves $C_t$ and $C_t'$ are already proved in \cite{ACL}. We state them here for the following reason.
Among the six curves considered in Proposition \ref{propthm1} and Proposition \ref{propthm2}, the curve $C_b$ is rational cuspidal and the curves $C_0'$ and $C_{b'}'$ are rational {\it nearly cuspidal} curves, that is all their singularities are irreducible except one which has 2 branches, see Remarks \ref{rkB1} and \ref{rkNB1}. Hence the claims above are compatible with fact that a rational cuspidal curve of even degree is either free or nearly free, see \cite{DStRIMS}, and
a rational nearly cuspidal curve $C$ of even degree is free, nearly free or
an $MPOG$-curve, see \cite{minTjurina,DStRIMS}. The curve $C_0$, which is not cuspidal, shows moreover that the hypothesis of being cuspidal is necessary in general for the first result to hold, even if there are exceptions, as shown by the curves $C_0'$ and  $C_{b'}'$ which are not cuspidal. The curves $C_1$ and $C_1'$ are both $MPOG$-curves, even though they are not rational.
 \end{rk}
 
 Next with discuss the curve arrangements involving some of the curves
 $C_t$ (resp. $C'_t$) and the high order contact lines $L_x:x=0$, $L_y:y=0$ and $L_z:z=0$.

  \begin{prop}
\label{propthm3}

\begin{enumerate}

\item The curves $C_{0} \cup L_z $ and $C_{1} \cup L_z$ are nearly free  with exponents $(5,6)$ and $\tau(C_{0} \cup L_z)= \tau(C_{1} \cup L_z)= 74$.

\item The curve $C_{b} \cup L_z$ is  free  with exponents $(4,6)$ and $\tau(C_{b} \cup L_z)=76$.

\item The curve $C_0 \cup C_{b}$ (resp. $C_0 \cup C_{b} \cup L_z$) is  free with exponents
$(9,10)$  and $\tau(C_0 \cup C_{b})=271$, (resp. exponents $(9,11)$) and 
$\tau(C_0 \cup C_{b}\cup L_z)=301$.

\end{enumerate}

 \end{prop}

   \begin{prop}
\label{propthm4}

\begin{enumerate}

 \item The curves   $C_0' \cup L_z$ and  $ C'_{b'} \cup L_z$ are free with exponents $(5,5)$ and $\tau(C'_{b'}\cup L_z)=75$.

 \item The curves $C_1' \cup L_z$   
 is nearly free with exponents $(5,6)$ and $ \tau(C_1' \cup L_z)=74$.

 \item The curve $C_0' \cup C'_{b'}$ (resp. $C_0' \cup C'_{b'}\cup L_z$) is   free with exponents $(9,10)$ and $\tau(C_0' \cup C'_{b'})=271$,
 (resp. exponents $(9,11)$ and $\tau(C_0' \cup C'_{b'}\cup L_z)=301$).
 
 \end{enumerate}

 \end{prop}

Note that the lines $L_x:x=0$ and $L_y:y=0$ satisfy 
$$|C_0 \cap L_x|=|C_0' \cap L_x|=2 \text{ and } |C_0 \cap L_y|=|C_0' \cap L_y|=3,$$
 hence they are also high order inflectional lines as for the curves $C_0$ and $C_0'$.  Recalling our constructions in \cite{DIPS}, it is natural to ask what happens if we add these lines to $C_0$ and to $C_0'$.
By adding the lines $L_x$ and $L_y$ to the curve $C_0$ we get new examples of maximal Tjurina curves of type $(d,r)=(2r-1,r)$ as defined above. The complete result is the following.
   \begin{prop}
\label{propthm5} 
 The curve  $C_0 \cup L_x$ is nearly free with exponents $(5,6)$, the curve $C_0 \cup L_y$ is  maximal Tjurina  of type $(11,6)$ and the curve $C_0 \cup L_x \cup L_y \cup L_z$ is  maximal Tjurina  of type $(13,7)$.
 
 The curve $C_0' \cup L_x$ is free with exponents $(5,5)$,  the curve $C_0' \cup L_y$ is nearly free with exponents $(5,6)$ and the curve $C_0' \cup L_x \cup L_y \cup L_z$ is nearly free with exponents $(6,7)$.
  \end{prop}

\subsection{The proofs of Propositions \ref{propthm1}, \ref{propthm2}, \ref{propthm3},  \ref{propthm4} and  \ref{propthm5}}

All these proofs involve the use of some computations done using a computer algebra software, for instance SINGULAR \cite{Sing}.
One can proceed in several ways, as follows.
\begin{enumerate}

 \item Use the defining equations given in \eqref{eqB2}, \eqref{eqB3}, \eqref{eqNB1} and \eqref{eqNB2} and list a minimal set of generators for the corresponding module of Jacobian syzygies. Then the values for the total Tjurina numbers can be obtained using Corollaries \ref{corCTC} and \ref{corCTC2}, or our new Theorem \ref{thm0}.
 
 \item Use the defining equations given in \eqref{eqB2}, \eqref{eqB3}, \eqref{eqNB1} and \eqref{eqNB2} to compute the minimal degree of a Jacobian syzygy in each case, compute the total Tjurina numbers as in Remarks \ref{rkB1} and \ref{rkNB1} and then use Corollaries \ref{corCTC} and \ref{corCTC2}, or  Theorem \ref{thm0} to decide whether the curve is free, or nearly free, or an MPOG-curve.
\end{enumerate}
One may use both approaches to check the results given by the
computer algebra software. 

In the pencils of curves $\al C_0 +\be C_b$ (resp. $\al C_0'+\be C_{b'}'$) there are 3 exceptional members (i.e., members more singular than the general member), namely $C_0,C_b$ and $10L_z$ (resp. $C_0'$, $C_{b'}'$. and $10L_z$).  Therefore the freeness of the curves $C_0 \cup C_{b}$ and  $C_0 \cup C_{b}\cup L_z$ (resp. $C_0' \cup C_{b'}'$ and  $C_0' \cup C_{b'}'\cup L_z$)
can be related to the results in \cite{Di+}, see also \cite{Mich,JV}.
In fact, the claim (3) in 
Propositions \ref{propthm3} and \ref{propthm4}  can be obtained from
\cite[Theorem 3.6]{Di+} stated here in Theorem \ref{thmE1}, as explained in Section 5, see in particular Corollary \ref{corE1} and Example \ref{exE1}.

 \begin{rk}
\label{rkE2a}
Consider the a 1-parameter version of the polynomial $g$ in
\eqref{eqB1}, namely
$$g(b)=p^2u-\frac{5}{3}ps+bs,$$
where $b \in \C$. If we apply the SINGULAR command $grobcov$ to the Jacobian ideal of the homogenization $f_0(b)$ of the polynomial $g(b)$, we get the following special values of the parameter $b$, and in each case we have determined the freeness type of the associated curve $C_0(b):f_0(b)=0$.
\begin{enumerate}

\item for $b=0$ or $9b^2-66b+25=0$ or $3b+10=0$ or $8b-5=0$, the  curve $C_0(b)$ is nearly free with exponents $(5,5)$ and $\tau(C_0(b))=60$.
 
\item for $45b^2-42b+125=0 $ or $19440b^4+20952b^3+69039b^2-54090b-92125=0$ or $48b+65=0$ or $3b+1=0$, the  curve $C_0(b)$ is MPOG with exponents $(5,5,6)$ and  $\tau(C_0(b))=59$.
\end{enumerate} 
It is interesting that only for the value $b=-\frac{1}{3}$ among all the above special values, the corresponding
function $g(b):\C^2 \to \C$ has no singularities. This can be determined by checking whether $1$ belongs to the Jacobian ideal $J_{g(b)}$ in $R=\C[x,y]$.
However, other special values for $b$ give rise to interesting curves as well.  Recall that a curve is said to be of type $(d,r,m)$ if it is an $m$-syzygy curve of degree $d$ with exponents $d_1=d_2= \ldots =d_m=r$, see \cite{EqDeg} for notation and more on such curves. For a curve of type $(d,r,m)$, we set
$$\Delta m=2r-d+3-m.$$
It is known that $\Delta m \geq 0$ and the equality holds exactly for Tjurina maximal curves, see \cite{maxTjurina}. For instance the value  $b'=\frac{5}{8}$ produces the following curves.
\begin{enumerate}

\item The curve $f_0(b')-x^{10}=0$ and the curve $f_0(b')-x^9z=0$ are curves of type $(10,9,11)$, hence  maximal Tjurina curves since $\Delta m=0$.
 
\item The curve $f_0(b')-x^4z^6=0$ is a curve of type $(10,7,6)$, with
$\Delta m=1$.

\item The curve $xyz(f_0(b')-x^{10})=0$ is a curve of type $(13,11,10)$, with
$\Delta m=2$.

\end{enumerate} 

Moreover, applying the SINGULAR command $grobcov$ to the Jacobian ideal of homogenization $f_0(b',c)$ of the polynomial 
$$g(b',c)=p^2u-\frac{5}{3}ps+b's-c,$$ 
with $b'=\frac{5}{8}$, we get  4 special values of the parameter $c$,
and two of them give rise to the following two free curves
$$p^2u-\frac{5}{3}psz^5+\frac{5}{8}sz^8+ \frac{25}{27}z^{10}=0$$
and
$$p^2u-\frac{5}{3}psz^5+\frac{5}{8}sz^8- \frac{163}{180}z^{10}=0,$$
both with exponents $(4,5)$.

\end{rk}

\begin{rk}
\label{rkB10}
The curve $C_0$ can be used to construct two counter-examples to
\cite[Conjecture 1.10]{DIPS} by producing a supersolvable curve which is not free as follows. 
For the first counter-example, consider the point $q=(0:1:0) \in C_0$ as above,
which is a point of multiplicity $6$ on the curve $C_0$. The line $L_z:z=0$ is an exceptional line passing through $q$ for the curve $C_0$, since
$$L_z \cap C_0= \{p,q\}$$
while a generic line $L$ through $q$ satisfies
$$|L \cap C_0|=\deg C_0 - \mult_q C_0+1=10-6+1=5.$$
The lines through $q$ different from $L_z$ have an equation $L_t:x=tz$
for $t \in \C$. To determine the points in $L_t \cap C_0$ distinct from $q$, one has to solve the following equation which is obtained from
\eqref{eqB2} by setting $x=t$ and $z=1$.
\begin{equation}
\label{eqB10}
t^6y^4+(4t^5+3t^4)y^3+(6t^4+\frac{19}{3}t^3+3t^2)y^2+(4t^3+\frac{11}{3}t^2+2t+1)y+t^2+\frac{1}{3}t-1=0.
\end{equation}
The discriminant of this equation in $y$ is
$$\Delta(t)=t^{20}c(t), \text{ where } c(t)=32768t^3-768t^2+1824t-243.$$
The exceptional line corresponding to $t=0$ is $L_x:x=0$, and the equations of the three exceptional lines corresponding to the three roots of $c(t)=0$ are just the 3 distinct factors of the cubic form
$$h(x,z)=32768x^3-768x^2z+1824xz^2-243z^3.$$
Consider now the curve
$$D_0:xzh(x,z)f_0=0.$$
The point $q$ is a {\it modular point} for the curve $D_0$, as any line $L'$ through $q$ distinct from an irreducible component of $D_0$ is a generic line for $D_0$ in the sense that
$$|L' \cap D_0|=\deg D_0 - \mult_q D_0+1=15-11+1=5.$$
A curve is  {\it supersolvable}  if it has a modular point, see \cite[Definition 1.9]{DIPS} and hence $D_0$ is supersolvable. A direct computation with SINGULAR shows that
$D_0$ is a 4-syzygy curve with exponents $(8,8,8,8)$. Therefore $D_0$ is a maximal Tjurina curve of type $(d,r)=(15,8)$ as in Remark \ref{rkNU}, and hence $\nu(D_0)=2$, a small freeness defect.

The second example of a supersolvable curve which is not free may be constructed using the same curve $C_0$ and the point $q'=(0:0:1)\notin C_0$.
The lines $L_x:x=0$ and $L_y:y=0$ are clearly exceptional for $C_0$,
and the remaining lines through $q'$ are given by the equation $L_t': y=tx$ with $t \in \C^*$. To determine the points in $L_t' \cap C_0$ distinct from $q'$, one has to solve the following equation which is obtained from
\eqref{eqB2} by setting $y=tx$ and $z=1$.
\begin{equation}
\label{eqB11}
t^4x^{10}+4t^3x^8+3t^3x^7+6t^2x^6+\frac{19}{3}t^2x^5+t(3t+4)x^4+\frac{11}{3}tx^3+(2t+1)x^2+\left(t+\frac{1}{3}\right)x-1=0.
\end{equation}
The discriminant of this equation in $x$ is
$\Delta'(t)=t^{39}c'(t)$, where
$$c'(t)= 282429536481t^5+276496482330144t^4+2414080421160192t^3+$$
$$+16059343010660352t^2+2540256075186176t+91534343012352.$$
The equations of the five exceptional lines corresponding to the five roots of $c'(t)=0$ are just the 5 distinct factors of the quintic form
$$h'(x,y)=x^5c'\left(\frac{y}{x}\right).$$
The curve
$$D_0':xyh'(x,y)f_0=0$$
has $q'$ as a modular point, and a direct computation with SINGULAR shows that
$D_0'$ is a 5-syzygy curve with exponents $(10,11,11,11,11)$ and $\nu(D'_0)=13$, a large freeness defect.

\end{rk}


\section{On higher degree Brian\c con-type polynomials}

 The authors in \cite{ACL} have constructed (several) Brian\c con-type polynomials $g_n$ for any degree $d_n=6n+4$, the polynomials $g$ and $g'$ in \eqref{eqB1} being the polynomials corresponding to $n=1$.
They are given by the formula
\begin{equation}
\label{eqHD1}
g_n=p^{2n}u+s\left ( \sum_{j=0}^na_jp^j+ \sum_{j=n+1}^{2n-1}p^j\right ),
\end{equation}
where $s,p,u$ are as in \eqref{eqB1} and the coefficients $a_0, \ldots, a_n$ take specified values in $\C$.
For $n=2$, there are 3 such polynomials constructed in \cite{ACL},
the simplest good choice (corresponding to Case 1 in \cite{ACL})  of the triple of coefficients being
\begin{equation}
\label{eqHD2}
a_0= -\frac{1}{5}, \  a_1= -\frac{3}{5}  \text{ and } a_2= -\frac{11}{5}.
\end{equation}
The corresponding polynomial, to be denoted by $g''$, has 2 atypical values, $0$ and 
$$b''=-(a_2-1)^2/4= -\frac{64}{25}.$$
The corresponding curves in $\PP^2$, denoted by $C_0''$, $C_{b''}''$ and $C_1''$ have similar properties to the curves $C_t$ and $C_t'$ discussed above. The equations for these curves are the following.

$$C_0'':f_0''=x^{10}y^6+6x^9y^5z^2+5x^8y^5z^3+15x^8y^4z^4+25x^7y^4z^5+20x^7y^3z^6+10x^6y^4z^6+$$
$$+50x^6y^3z^7+15x^6y^2z^8
+\frac{184}{5}x^5y^3z^8+50x^5y^2z^9+10x^4y^3z^9+6x^5yz^{10}+$$
$$+\frac{252}{5}x^4y^2z^{10}+25x^4yz^{11}+22x^3y^2z^{11}+x^4z^{12}+\frac{152}{5}x^3yz^{12}+5x^2y^2z^{12}+$$
$$+5x^3z^{13}+14x^2yz^{13}+\frac{34}{5}x^2z^{14}+4xyz^{14}+2xz^{15}+yz^{15}-z^{16}=0,$$
$$C_1'':f_1''=f_0''-z^{16}=0 \text{ and } C_{b''}:f_{b''}''=f_0+\frac{64}{25}z^{16}=0.$$
It was shown in \cite{ACL} that $g(C_0'')=g(C_{b''}'')=0$, $g(C_1'')=2$
and all the curves $C_t''$ have two points on the line at infinity $L$.
The interested reader can check, using one of the approaches $(1)$ or $(2)$ explained above at the end of Section 3, the following result. 
Moreover, the claim (3) in 
Proposition \ref{propthm6}  can be obtained from
\cite[Theorem 3.6]{Di+} as explained below in Corollary \ref{corE1} and Example \ref{exE1}.

   \begin{prop}
\label{propthm6}

\begin{enumerate}

 \item The curves $C_0''$ and $C_1''$ are MPOG-curves with exponents $(8,8,9)$ and $\tau(C_0'')=\tau(C_1'')=167$. The curve $C_{b''}''$ is free with exponents
$(7,8)$ and $\tau(C_{b''}'')=169$. 

\item The curve $C_{b''}'' \cup L$  is free with exponents $(7,9)$ and $\tau(C_{b''}'' \cup L)=193$ and the curves
$C_{0}'' \cup L$ and $C_{1}'' \cup L$ are nearly free with exponents $(8,9)$ and $\tau(C_{0}'' \cup L)=\tau(C_{1}'' \cup L)=191$.

\item The curve $C_0'' \cup C_{b''}'' $ is free with exponents $(15,16)$ and $\tau(C_0'' \cup C_{b''}'')=721$ and the curve $C_0'' \cup C_{b''}'' \cup L$ is free with exponents $(15,17)$ and $\tau(C_0'' \cup C_{b''}'')=769$.

\end{enumerate}
 \end{prop}
 
 For each $n \geq 2$, the authors  in \cite{ACL} have shown the existence of a specific choice of the coefficients in \eqref{eqHD1}, corresponding to the Case $n+1$ in \cite{ACL}. Explicit formulas for the polynomial $g_n'$ obtained from this choice are given in \cite{ACL} for the cases $2 \leq n \leq 4$,
which involve rational numbers with large denominators, see Remark \ref{rk7}.

 For this polynomial $g_n'$, the projective closure $C_{n,0}'$ of the affine fiber $g_n^{' -1}(0)$ is a rational curve, having 3 branches at infinity by \cite[Figure 7]{ACL}. Since $C_{n,0}'$ has  2 points at infinity, it follows that  $C_{n,0}'$ is a nearly cuspidal rational curve, as defined in Remark \ref{rkG} above. Using the proof of \cite[Theorem 1.3]{minTjurina} and the statement of \cite[Theorem 3.11]{minTjurina}, we obtained the following.

   \begin{thm}
\label{thm7}
For any $n \geq 2$, the curve $C_{n,0}'$ coming from a special fiber of the Brian\c con-type polynomial $g_n'$ is either free, or nearly free or an MPOG-curve.
 \end{thm}
   \begin{rk}
\label{rk7} 
 A direct computation with SINGULAR shows that $C_{2,0}'$ is a free curve of degree 16 with exponents $(7,8)$, $C_{3,0}'$ is a free curve of degree 22 with exponents $(7,14)$ and $C_{4,0}'$ is a free curve of degree 28 with exponents $(7,20)$.
 The choice of the coefficients $a_j$'s, that gives rise to the polynomial $g_n'$ are listed in \cite{ACL} as follows.
 For $n=2$ the coefficients are
 $$a_0= -\frac{1}{125}, \  a_1= \frac{17}{125}  \text{ and } a_2= -\frac{91}{125},$$
  for $n=3$ the coefficients are
$$ a_0= \frac{1}{2401}, \  a_1= -\frac{31}{2401}, \  a_2= \frac{353}{2401}   \text{ and }  a_3=-\frac{1695}{2401}$$
and the coefficients for the case $n=4$ 
are the following
$$
 a_0= -\frac{1}{N}, \  a_1= \frac{49}{N}, \  a_2= -\frac{951}{N} , \    a_3=\frac{9049}{N} \text{ and }  a_4=-\frac{40951}{N},
$$
where $N=59049$. We do not know whether the curves $C_{n,0}'$ are free for any $n \geq 5$.
 
 \end{rk} 
 
 \section{Eigenscheme ideals, freeness properties and pencils} 
 
In this section, we prove first Theorem \ref{thmE2}.
We start by showing  that $r={\rm mdr}(f)$. Let $d_1\leq d_2$ be first two  exponents of $C$ as in Introduction. Then $d_1 \leq r$. Since the eigenscheme $\Gamma _{\theta}$ is 0-dimensional, it follows that the derivation $\theta$ is primitive, namely its coefficients $a,b,c$ do not have a common factor. This implies that either $d_1=r$ or that $d_2 \leq r$. In the latter case, we have
 $$d_1+d_2 \leq r+r \leq d-1.$$
 But for any plane curve $C$ the reversed inequality $d_1+d_2 \geq d-1$ holds, with equality precisely when $C$ is a free curve. Hence in both cases $d_1=r$, and our claim is proved.

 The proof of  \cite[Theorem 2.5]{Di+} shows if $f \in I_{\theta}$, then there is a derivation $\mu \in D(f)$ such that the determinant of the $3 \times 3$ matrix $M(E,\theta, \mu)$ having as first row the coefficients
 of the Euler derivation $E$, as second row the coefficients
 of the derivation $\theta$ and as third row the coefficients
 of the  derivation $\mu$ is equal to $f$. Using the decomposition
 \eqref{eqE1} we may assume that $\mu \in D_0(f)$.
 
 Now replace the assumption $f \in I_{\theta}$ by the assumption
 $(K_{\theta,f})_1\ne 0$. This means that there exists a nonzero linear form $\ell \in S_1$ such that $\ell f \in I_{\theta}$. Then exactly the same argument as in the proof of  \cite[Theorem 2.5]{Di+} show the existence of a derivation
 $\eta \in D_0(f)$ such that
 \begin{equation}
\label{eqE2}
 \det M(E, \theta, \eta)= \ell f,
 \end{equation}
 where $M(E, \theta, \eta)$ is the $3 \times 3$ matrix constructed using the coefficients of $E$, $\theta$ and $\eta$ as above. Since $\ell f \ne 0$, it follows that $\eta$ is not a multiple of $\theta$ and hence $d_2 \leq \deg \eta=d-r$, where $d_2$ is the second exponent of $C$ as in Introduction. In this way we get $r+d_2=d_1+d_2 \leq d$.
 By assumption $K_{\theta,f}$ is a proper ideal, hence $C$ is not a free curve, and this implies $r+d_2=d_1+d_2 >d-1$. It follows that $d_1+d_2=d$ and hence $C$ is a plus-one generated curve by \cite[Theorem 2.3]{minTjurina}. Conversely, if $C$ is a plus-one generated curve and we take $\eta$ the derivation corresponding to a generator of $D_0(f)$ of degree $d_2=d-r$, then the equality \eqref{eqE2} holds for some linear form $\ell$, see Remark \ref{rkE1} below if necessary.
This yields $\ell \in (K_{\theta,f})_1$ and hence $(K_{\theta,f})_1\ne 0$.
 Since $C$ is a plus-one generated curve, $C$ is a 3-syzygy curve. If we denote the third generator in a minimal set of such generators by $\rho_3$ as in Introduction, and by $\eta_3$ the corresponding derivation in $D_0(f)$, then the second generator of $K_{\theta,f}$ is exactly
 $$h= \det M(E, \theta, \eta_3)/f$$
 and all the claims in Theorem \ref{thmE2} are now proven.
 \endproof 
   \begin{rk}
\label{rkE1} 
The ideal $K_{\theta,f}$ defined in Theorem \ref{thmE2} clearly coincide
by the remarks at the beginning of the above proof with the Bourbaki ideal $B(C,\rho)$, where $\rho \in D_0(f)$ is the syzygy associated to the derivation $\theta$, see for instance Section 3 in \cite{minTjurina}
or \cite{JNS} for more on such Bourbaki ideals.
In particular, \cite[Proposition 3.1]{minTjurina} is related to the final part of the above proof.
 \end{rk}  
 Theorem \ref{thmE1} is used in \cite{Di+} to construct free curve arrangements obtained by taking the union of some members in a pencil of plane curves. In this section we recall their main result, see \cite[Theorem 3.6]{Di+}, in a special situation, which is enough for our purpose.
 
 Let $C:f=0$ be a reduced curve in $\PP^2$ of degree $d$ and consider the pencil of curves 
 $$C_u: f_u=\al f(x,y,z)+ \be z^d=0,$$
 where $u=(\al:\be) \in \PP^1$. It is clear that the derivation
 $$\delta=\delta_{f,z}=f_{y}\partial_x- f_{x}\partial_y$$
 belongs to $D_0(f_u)$ and in fact to $D_0(F_k)$, where
 $F_k$ is the product of $k \geq 1$ distinct reduced member of the pencil $C_u$. With this notation, we set $F'_k=zF_k$ and note that 
 $\delta \in D_0(F'_k)$ as well. Assume from now on that
 $$ (A) \  \   \text { the partial derivatives } f_x \text { and } f_y \text { have no common factor}.$$
 Then $\delta$ is a primitive derivation, i.e. it is not the multiple of a strictly lower degree derivation in $S$. It follows that either $\delta$
 is a minimal degree derivation for $F_k$, or one has
 $$d_1 \leq d_2 \leq \deg \delta =d-1,$$
 where $d_1$ and $d_2$ denote the first two exponents of $F_k$.
 When $k \geq 2$, the second case is impossible, since it leads to the contradiction
 $$d_1+d_2 \leq 2(d-1) <\deg F_k -1.$$
 Hence, when $k\geq 2$ we have that $\delta$ is a minimal degree derivation for $F_k$ and $F_k'$, with 
 $$\deg \delta =d-1={\rm mdr}(F_k)={\rm mdr}(F'_k).$$
 Consider now the ideal $I_{\delta}$ and note that the corresponding
 eigenscheme $\Gamma_{\delta}$ is 0-dimensional when our assumption $(A)$ holds and in addition $C$ has not the line at infinity $L_z:z=0$ as an irreducible component.
 With this notation, we have the following consequence of Theorem \ref{thmE1}, a special case of \cite[Theorem 3.6]{Di+}.
 
\begin{cor}
\label{corE1} If the assumption $(A)$ holds and 
the line at infinity $L_z$ is not an irreducible component for the curve $C$,
then the curve arrangement $\CC_k:F_k=0$ (resp. $\CC_k':F'_k=0$) for $k \geq 2$ is free
 if and only if
$F_k \in I_{\delta}$ (resp. $F'_k \in I_{\delta}$). If this is the case, then the  exponents of $\CC_k$ (resp. $\CC'_k$) are $(d-1, (k-1)d)$ (resp.  $(d-1, (k-1)d+1)$.
  \end{cor}

  \begin{ex}
\label{exE1}
With the above notation,  consider the pencil $C_u$, where $f=f_0$ from \eqref{eqB2}. It is clear that both assumptions in Corollary \ref{corE1} are
satisfied.
Using SINGULAR it is easy to check that $F_2=f_0f_b \in I_{\delta}$ and respectively 
$F_2'=zf_0f_b \in I_{\delta}$.
This yields the claim (3) in Proposition \ref{propthm3}. 
The proof of the claim (3) in Proposition \ref{propthm4} is completely similar.
Moreover, one can check using SINGULAR that
 \begin{equation}
\label{eqE20}
f_0^2 \in I_{\delta}, \ z^{10}f_0 \in I_{\delta} \text{ and } z^{20} \in I_{\delta}
 \end{equation}
and also
$${f'}_0^2 \in I_{\delta'}, \ z^{10}f'_0 \in I_{\delta'} \text{ and } z^{20} \in I_{\delta'},$$
where $\delta'=f'_{0y} \partial_x- f'_{0x}\partial_y $, with $f_0'$ defined in \eqref{eqNB1}.
It would be nice to find a theoretical explanation for the relations above, but this does not seem to be easy. We note here only that the fact that $g,g':\C^2 \to \C$ have no singularities can be deduced from the following relations
$$  ( 4x^4y +4x^3 +x^2)g_x-(  6x^3y^2 +8x^2y +2x -1)g_y=1$$
and 
$$-(48x^6y^2 +96x^5y +72x^4y +48x^4 +\frac{88}{3}x^3 +15x^2)g'_x+  $$
$$+(72x^5y^3 +168x^4y^2 +90x^3y^2 +120x^3y +44x^2y +24x^2 -\frac{10}{3}x +1)g'_y=1,                                   $$
see also \cite{ACL0} for the first relation. By homogeneization,
 the first relation implies that $z^{14} \in I_{\delta}$ and the second relation implies that $z^{17} \in I_{\delta'}$.

It follows by  Corollary \ref{corE1} and using \eqref{eqE20} that any union  
$$\CC_k=C_{t_1} \cup \ldots \cup C_{t_k}$$
of $k\geq 2 $ distinct members of the pencil $C_t: f_0 + t z^{10}=0$, with 
$t _j \in \C$ gives rise to a free curve with exponents $(9,10k-10)$.
Similarly, the union $\CC_k \cup L_z$ is free with exponents $(9, 10k-9)$.
The same claims, with the same proof, hold for the pencil $C_t': f_0' + t z^{10}=0$.

\end{ex}
  \begin{rk}
\label{rkE2}
The curve $C_b \cup L_z$ from Proposition \ref{propthm3} is free with exponents $(4,6)$. One cannot apply Corollary \ref{corE1} to this curve, since in this case the number of reduced member of the pencil used in the arrangement is $k=1$. On the other hand,  \cite[Theorem 3.6]{Di+} applies to this curve, and the equivalent claims in the conclusion of that result both fail in this case. This shows that Corollary \ref{corE1} is a new version of a special case of  \cite[Theorem 3.6]{Di+}.
\end{rk}
Finally we show that there is an analog of Corollary \ref{corE1} in the case of MPOG-curve arrangements. Theorem \ref{thmE2} yields the following result.

\begin{cor}
\label{corE2} If the assumption $(A)$ holds and 
the line at infinity $L_z$ is not an irreducible component for the curve $C$ of degree $d$,
then the curve arrangement $\CC_k:F_k=0$ (resp. $\CC_k':F'_k=0$) for any $k \geq 2$ is MPOG
 if and only if the ideal $K_{\delta,F_k}$ (resp. $K_{\delta,F'_k}$)
has two generators, a linear form $\ell \in S_1$ and a quadratic form $q \in S_2$.
If this is the case, then the  exponents of $\CC_k$ (resp. $\CC'_k$) are $(d-1, (k-1)d+1, (k-1)d+2)$ (resp.  $(d-1, (k-1)d+2, (k-1)d+3)$).
  \end{cor} 
  The following example shows that the conditions in Corollary \ref{corE2}
  can be checked in practice.

With the above notation,  consider the pencil $C_u$, where $$f=f_0\left(\frac{5}{8}\right)=x^6y^4+4x^5y^3z^2+3x^4y^3z^3+6x^4y^2z^4+19/3x^3y^2z^5+4x^3yz^6+$$
$$+3x^2y^2z^6+\frac{11}{3}x^2yz^7+x^2z^8+\frac{71}{24}xyz^8+\frac{1}{3}xz^9+yz^9-\frac{1}{24}z^{10},$$ as in  Remark \eqref{rkE2a} above. Hence the curve 
$$C_0\left(\frac{5}{8}\right): f_0\left(\frac{5}{8}\right)=0$$
is nearly free with exponents $(5,5)$.
It is clear that both assumptions in Corollary \ref{corE2} are
satisfied for this curve.
Using SINGULAR it is easy to check that
$$K_{\delta, f^2}=K_{\delta, z^{10}f}=K_{\delta, z^{20}}=(\ell_1, \ell_2\ell_3),$$
where 
$$\ell_1=793202173x+3698829360y-1039006968z$$
and
$$\ell_2\ell_3=
106288200y^2-23416645yz-11726872z^2.$$
The ideal $K=(\ell_1, \ell_2\ell_3)$ is a complete intersection, defining a subscheme of $\PP^2$ consisting of two simple points, namely
$$p=V(\ell_1, \ell_2) \text{ and } q=V(\ell_1,\ell_3).$$
A polynomial $h \in S$ is in the ideal $K$ if and only if $h(p)=h(q)=0$.
Using this remark, it is easy to see that neither $f$ nor $z$, nor any reduced curve in the pencil $C_t(\frac{5}{8}): f + t z^{10}=0$ is in the ideal $K$. It follows that 
$$K_{\delta, f^p}=K_{\delta, z^{q}f^r}=K_{\delta, z^{s}}=K$$
for any $p \geq 2$, $q\geq 10$, $r \geq 1$ and $s \geq 20$.
This implies that for any $k \geq 2$ one has
$$K_{\delta, F_k}=K_{\delta, F'_k}=K.$$
Applying Corollary \ref{corE2}, this implies the following.

\begin{prop}
\label{propthmE3} 
Any union  
$$\CC_k=C_{t_1} \left(\frac{5}{8}\right)\cup \ldots \cup C_{t_k}\left(\frac{5}{8}\right)$$
of $k\geq 2 $ distinct reduced members of the pencil $C_t(\frac{5}{8}): f_0(\frac{5}{8}) + t z^{10}=0$, with 
$t _j \in \C$ gives rise to an MPOG curve with exponents $(9,10k-9,10k-8)$.
Similarly, the union $\CC'_k=\CC_k \cup L_z$ of $\CC_k$ with the line at infinity $L_z:z=0$ is an MPOG curve with exponents $(9, 10k-8,10k-7)$.
\end{prop}
As an explicit example for the curve arrangement $\CC_k$ (resp. $\CC'_k$ in Proposition
\ref{propthmE3} one may take
$$\CC_k: f_0\left(\frac{5}{8}\right)^k+z^{10k}=0$$
and respectivelly
$$\CC'_k: z\left(f_0\left(\frac{5}{8}\right)^k+z^{10k}\right)=0$$
for any $k \geq 2$.

\end{document}